\documentclass[12pt]{article}

\usepackage{array}
\usepackage{amsmath}
\usepackage{amssymb}
\usepackage{latexsym}
\usepackage{graphicx}
\usepackage[usestackEOL]{stackengine}

\usepackage{mathtools}
\usepackage{tikz-cd}
\usepackage{multicol}

\usetikzlibrary{positioning}
\usetikzlibrary{arrows}

\usepackage{rotating} 

\usepackage{xcolor}

\usepackage{url}

\usepackage{amsthm}
\newtheoremstyle{tssf}{10pt}{10pt}{\sffamily}{}{\sffamily}{:}{.5em}{}
\theoremstyle{tssf}
\newtheorem{theorem}{Theorem}[section] 
\newtheorem{proposition}[theorem]{Proposition}

\newtheorem{definition}[theorem]{Definition} 
\newtheorem{remark}[theorem]{Remark}
\newtheorem{example}[theorem]{Example}

\usepackage{xspace}

\usepackage{mathrsfs}

\newcommand{\VBgpds}{$\mathscr{V}\!\mathscr{B}$--groupoids\xspace}

\newcommand{\dvb}{double vector bundle\xspace}
\newcommand{\dvbs}{double vector bundles\xspace}
\newcommand{\vb}{vector bundle\xspace}
\newcommand{\vbs}{vector bundles\xspace}

\newcommand{\warp}{w}           
  
\newcommand{\dual}{(*)}

\newcommand{\add}[1]{\mathbin{\lower 5pt%
    \hbox{${\stackrel{\textstyle +}{\scriptscriptstyle #1}}$}}}
\newcommand{\sub}[1]{\mathbin{\lower 5pt%
    \hbox{${\stackrel{\textstyle -}{\scriptscriptstyle #1}}$}}}
\newcommand{\sm}[1]{\mathbin{\lower 5pt
    \hbox{${\stackrel{\textstyle \hdot}{\scriptscriptstyle #1}}$}}}
\newcommand{\by}[1]{\mathbin{\lower 7pt
    \hbox{${\stackrel{\cdot}{\scriptscriptstyle #1}}$}}}

\renewcommand{\Bar}[1]{\overline{#1}}
\renewcommand{\Tilde}[1]{\widetilde{#1}}

\newcommand{\vfk}[1]{\mathscr{X}(#1)}

\usepackage{stmaryrd}
\usepackage{wasysym}

\newcommand{\D}{\mathop{\raise0.1ex\hbox{$\mathfrak{D}$}}}

\newcommand{\kduer}{%
\mathbin{\raisebox{3pt}{\varhexstar}\kern-3.70pt{\rule{0.15pt}{4pt}}}\,}
 \newcommand{\duer}{\mathbin{\raisebox{3pt}{\varhexstar}\kern-3.85pt{\rule{0.24pt}{5pt}}}\,}

\newcommand{\act}{\mathbin{\hbox{$<\kern-.4em\mapstochar\kern.4em$}}}
\newcommand{\ract}{\mathbin{\hbox{$\mapstochar\kern-.3em>$}}}

\newcommand{\tilo}{\widetilde 0}           

\newcommand{\llangle}{\langle\!\langle}
\newcommand{\rrangle}{\rangle\!\rangle}

\newcommand{\lpair}{\angle}
\newcommand{\rpair}{\reflectbox{$\angle$}} 

\newcommand{\thbr}[2]%
{\rule[-1pt]{1pt}{10pt}\hspace{2pt} #1,\, #2\hspace{1pt}\rule[-1pt]{1pt}{10pt}}

\newcommand{\newpair}[2]%
{\talloblong\hspace{1pt} #1,\, #2\hspace{0.5pt}\talloblong}

\newcommand{\R}{\mathbb{R}}

\newcommand{\id}{\text{\upshape id}}

\newcommand{\Ga}{\Gamma}

\renewcommand{\phi}{\varphi}
\newcommand{\ph}{\varphi}

\newcommand{\goX}{\mathfrak{X}}
\newcommand{\goF}{\mathfrak{F}}

\newcommand{\hlift}{H}
\newcommand{\co}{:}

\newcommand{\st}{\ \vert\ }

\DeclareMathOperator*{\equals}{=} 

\usepackage{subfig}
\usepackage{caption}
\newcommand{\RM}[2]{\mx^{#1}_{#2}}

\usepackage{mathrsfs}
\newcommand{\mscr}[1]{\mathscr{#1}}

\newcommand{\mA}{\mscr{A}}
\newcommand{\mB}{\mscr{B}}

\usepackage{MnSymbol}

\newcommand{\mx}{R} 

\newcommand{\nsp}{nonstandard pairing\xspace}

\newcommand{\lhangle}{\,\rule[-0.1ex]{0.4ex}{0.7\baselineskip}\,\,}
\newcommand{\rhangle}{\,\,\rule[-0.1ex]{0.4ex}{0.7\baselineskip}\,}

\newcommand{\sol}{\bullet}

\begin{document}

\sffamily

\title{\textbf{Warps and duality\\ for \dvbs}}

\author{Magdalini K. Flari\footnote{Max Planck Institute for Mathematics, 
Vivatsgasse 7, 53111 Bonn, Germany, and Department of Mathematics, 
Aristotle University of Thessaloniki, Thessaloniki, 54124, Greece,
  \text{flari@math.auth.gr}}
  and
Kirill Mackenzie\footnote{School of Mathematics and Statistics, 
  University of Sheffield, Sheffield, S3 7RH, United Kingdom,
\text{K.Mackenzie@sheffield.ac.uk}}}

\date{}
	
\maketitle
	
\begin{abstract}
  A linear section of a \dvb is a parallel pair of sections
  which form a \vb morphism; examples include the complete
  lifts of vector fields to tangent bundles and the horizontal
  lifts arising from a connection in a \vb. A grid in a \dvb
  consists of two linear sections, one in each structure, and
  thus provides two paths from the base manifold to the top
  space; the warp of the grid measures the lack of commutativity
  of the two paths.

  For a \dvb $D$, a linear section induces a linear map on the relevant
  dual of $D$ and thence a linear section of the iterated dual; we call
  this the corresponding \emph{squarecap section}. This notion makes precise
  the relationships between various standard $1$-forms and vector fields
  that are not dual in the usual sense. A grid in $D$ induces squarecap
  sections in the two iterated duals. These iterated duals are themselves
  in duality and we obtain an expression for the warp of a grid in terms
  of a pairing. This result extends the relationships between various standard
  formulas, and clarifies the iterated dualization involved in the concept of
  double Lie algebroid. 

  In the first part of the paper we review the nonstandard pairing
  of the two duals of a \dvb, showing in particular that it arises
  as a direct generalization of the Legendre-type antisymplectomorphism
  $T^*(A^*)\to T^*(A)$ for any \vb $A$ of Mackenzie and Xu. 
\end{abstract}

\vfil

\newpage 
 
\section{Introduction}
\label{sect:intro}

This paper is concerned with two aspects of the duality theory
of \dvbs. Firstly, we show that the non-degenerate pairing of
the two duals of a \dvb \cite{Mackenzie:1999}
emerges naturally from the canonical antisymplectomorphism
$\mx\co T^*(A^*) \to T^*(A)$, for $A$ a \vb, 
introduced by Ping Xu and Mackenzie \cite{MackenzieX:1994}.
This map, which may be regarded as a Legendre transformation,
is crucial to work on Lie bialgebroids and related concepts in
Poisson geometry. The non-degenerate pairing of the duals is,
likewise, crucial to the concept of double Lie algebroid 
as introduced in \cite{Mackenzie:2011}. 

Secondly we show that by using the duality properties of \dvbs, and
in particular the pairing of the duals, the warp of a grid in a \dvb
may be expressed in terms of the dual structures. A grid in a \dvb
consists of sections of each of the four \vb structures, such that
each pair of parallel sections forms a \vb morphism,  
see \cite[Definition 1]{FlariM:2019}. 
The warp of a grid measures the lack of commutativity of the four 
sections. 
Well-known geometric objects can be expressed as warps, for
example, the bracket of two vector fields is a warp, and, given
a connection in a vector bundle, the covariant derivative of a
section along a vector field is a warp, see \cite[Example 1]{FlariM:2019}.

The concept of the dual \dvb of a \dvb is a particular case of the duality
of \VBgpds, introduced by Pradines \cite{Pradines:1988}, in what was 
the first construction of the dual of a double structure. The non-degenerate
pairing of the two duals, which we call the \emph{\nsp}, was originally defined
by a simple explicit formula (\ref{eq:nspAB}); we demonstrate in
\S\ref{sect:mx} that it is a direct generalization of the construction
of $\mx$.

The paper supposes a familiarity with the notion of a \dvb, and
of its duals, such as in \cite[\S9.2]{Mackenzie:GT}, but a
detailed knowledge of the duality theory 
is not needed. We give a brief review here.

Start with a \dvb $D$ as shown in the first diagram in Figure~\ref{fig:Dds}.
Taking the dual (in the usual sense) of $D$ with respect to the vector bundle structure
with base $A$ defines another \dvb, denoted $D\duer A$, as in the second diagram in
Figure~\ref{fig:Dds}. Here $C$ is the core of $D$, defined as the set of elements of
$D$ that project to zeros in both $A$ and $B$; the two additions coincide on $C$ and
give it the structure of a \vb on base $M$. The core of $D\duer A$ is $B^*$.

In a similar way we can define the \dvb $D\duer B$, with core $A^*$. 
The fourth diagram in Figure~\ref{fig:Dds} is an iterated dual; we will need it later. 
\begin{figure}[h]
$$
\begin{tikzcd}[row sep=1cm, column sep = 1cm]
D \arrow[r, ""] \arrow[d,""] & B \arrow[d,""] \\
A\arrow[r,""] &M,
\end{tikzcd}
\hspace*{8mm}
\begin{tikzcd}[row sep=1cm, column sep = 1cm]
D\duer A \arrow[r, ""] \arrow[d,""] &C^* \arrow[d,""] \\
A\arrow[r,""] &M,
\end{tikzcd}
\hspace*{8mm}
\begin{tikzcd}[row sep=1cm, column sep = 1cm]
D\duer B \arrow[r, ""] \arrow[d,""] & B \arrow[d,""] \\
C^* \arrow[r,""] &M.
\end{tikzcd}
\quad
\begin{tikzcd}[ampersand replacement=\&, row sep = 1cm, column sep = 1cm]
D\duer B\duer C^* \arrow[r, ""]
\arrow[d,"",swap]
\& A \arrow[d,""]
\\
C^*\arrow[r,"",swap] \& M.
\end{tikzcd}
$$
\caption{\label{fig:Dds}}
\end{figure}

There exists a non-degenerate pairing between $D\duer A$ and $D\duer B$ over
$C^*$ \cite{Mackenzie:1999}. This pairing is natural up to sign and in both
\cite[9.2.2]{Mackenzie:GT} and \cite[2.2]{Mackenzie:2011} it is defined by
\begin{equation}
  \label{eq:nspAB}
  \lhangle\Phi,\Psi\rhangle
  = \langle \Phi,d\rangle_A-\langle \Psi, d\rangle_B,
\end{equation}
where $\Phi\in D\duer A$, and $\Psi\in D\duer B$, with outlines
$(\Phi;a,\kappa;m)$ and $(\Psi;\kappa,b;m)$, and where
$d$ is any element of $D$ with outline $(d;a,b;m)$. 
We refer to the choice of signs on the right hand side as the
\emph{`AB convention'}.

We refer to this pairing, either in this convention or with
the opposite sign convention, as the \emph{\nsp}. 

Associated with the \nsp in \cite{Mackenzie:GT,Mackenzie:2011}, are
two maps $Z_A$ and $Z_B$; we consider these in Remark~\ref{rmk:ZZ}. 

In \S\ref{sect:mx} we show how the \nsp arises naturally from
the $\mx$ map. 

The concepts of warp and grid are recalled in 
\S\ref{sect:warpsandgrids}.
Only the basic definitions from \cite{FlariM:2019} are needed
in this paper. 

\subsection*{Acknowledgements}

The material in sections \ref{sect:mx} and \ref{sect:dvbsd} is based on
seminars given by Mackenzie at the University of Manchester, the IHP and
the MPIM in late 2018 and he is grateful to the organizers of these seminars
for these opportunities. The authors are very grateful to Fani Petalidou
for a very careful reading of the paper and for catching many imperfections.

\section{The canonical antisymplectomorphism $T^*(A^*)\to T^*(A)$}
\label{sect:mx}

In \cite{MackenzieX:1994} Ping Xu and Mackenzie established the following result.

\begin{theorem}
\label{thm:MX94}
Given a vector bundle $A$ on base $M$, 
there is a diffeomorphism
$$
\mx\co T^*(A^*)\to T^*(A)
$$
which is an antisymplectomorphism of the canonical cotangent structures. 
\end{theorem}

Further properties of $\mx$, which provide a characterization, are described below. 
The case $A = TM$ is due to Tulczyjew \cite{Tulczyjew:1977ihp}. 

A formulation for the cotangent bundles of super vector bundles was given by Roytenberg
\cite{Roytenberg:thesis}; see also \cite{KhudaverdyanV:2018}. 

We repeat the proof of \ref{thm:MX94} here since it is needed in \S\ref{sect:dvbsd}. 
The key is the following result of
Tulczyjew \cite{Tulczyjew:1977ihp,Tulczyjew:1977}. 

\begin{proposition}
  \label{prop:Tulc3.1}
  Let $S$ be a submanifold of a manifold $Q$ and let $F\co S\to\R$
  be smooth. Then 
  $$
  W := \{\phi\in T^*_sQ \st s\in S,\ \langle\phi,X\rangle = \langle dF,X\rangle
  \text{ for all } X\in T_sS\} 
  $$
  is Lagrangian in $T^*Q$. 
\end{proposition}

{\sc Proof of Theorem \ref{thm:MX94}:} In $Q:= A^*\times A$ define
$S$ to be $A^*\times_M A$ and $F\co S\to\R$ to be the standard pairing.

Consider pairs $(\goF,\Phi)\in T^*(A^*)\times T^*(A)$
that project to the same $m\in M$; say $\goF\in T^*_\psi(A^*)$ and
$\Phi\in T^*_X(A)$. Then $(\goF,\Phi)$ will lie in $W$ if for all
$(\goX,\xi)\in T_{(\psi,X)}(S) = T_\psi(A^*)\times_{TM} T_X(A)$, we
have 
\begin{equation}
  \label{eq:cond-v1}
\langle(\goF,\Phi),(\goX,\xi)\rangle =
\langle dF,(\goX,\xi)\rangle. 
\end{equation}
The right-hand side is the tangent pairing of $T(A^*)$ and
$T(A)$ over $TM$, denoted $\llangle\goX,\xi\rrangle$ in \cite{MackenzieX:1994}.
Expanding out the left-hand side, (\ref{eq:cond-v1}) becomes
\begin{equation}
\label{eq:cond-v2}
\langle\goF,\goX\rangle + \langle\Phi,\xi\rangle =
\llangle\goX,\xi\rrangle.
\end{equation}
That is, $W$ consists of those $(\goF,\Phi)$, projecting to the same $m\in M$,
such that (\ref{eq:cond-v2}) holds for all $(\goX,\xi)$ with $T(q)(\goX) = T(q_*)(\xi)$,
where $q$ and $q_*$ are the bundle projections. 

Note that $\goF$, $\goX$, $\Phi$ and $\xi$ in (\ref{eq:cond-v2}) all
project to the same $m\in M$. To show that, given $\goF$, there is a
unique $\Phi$ such that (\ref{eq:cond-v2}) holds, it is therefore
sufficient to consider the case $A = M\times V$. 

Write (adapting the notation for projections)
\begin{gather}
\goF = (\chi,\psi,Y)\in T^*M\times V^*\times V, \qquad
\Phi = (\omega, X, \ph)\in T^*M\times V\times V^*.\\
\goX = (x,\psi,\theta)\in TM\times V^*\times V^*,\qquad
\xi  = (x,X,Z)\in TM\times V\times V.
\end{gather}
Then
$$
\langle{\goF},{\goX}\rangle + \langle{\Phi},{\xi}\rangle
 = \langle\chi, x\rangle + \langle\theta, Y\rangle
 + \langle\omega,x\rangle + \langle\ph,Z\rangle
$$
and
$$
\llangle\goX,\xi\rrangle = \langle\psi,Z\rangle + \langle\theta, X\rangle.
$$
The variables $\theta$ and $Z$ do not take part in the compatibility
condition. We can therefore set $\theta = 0$ and $Z = 0$.
So $\omega = -\chi$, and we now have
$$
\langle\theta, Y\rangle + \langle\ph,Z\rangle = \langle\psi,Z\rangle + \langle\theta, X\rangle.
$$
Now setting $\theta = 0$ we get $\ph = \psi$.  Lastly, $X = Y$.

So $\Phi = (-\chi,Y,\psi)$ is determined by $\goF$. Reversing the
argument, $\goF$ is determined by $\Phi$. So $\mx$ is the graph of
a diffeomorphism and, locally, $R$ is $(\chi,\psi,Y)\mapsto (-\chi,Y,\psi)$.

Since $W$ is Lagrangian, $\mx$ is an antisymplectomorphism. Further
\cite[9.5.2]{Mackenzie:GT}, the Liouville 1-forms are related by
$\mx^*\lambda_A + \lambda_{A^*} = dP$, where $P\co T^*(A^*)\to\R$
sends $\goF$ to $\langle c_*(\goF), c(\mx(\goF)\rangle$; here
$c$ and $c_*$ are the standard projections for the cotangent
bundles $T^*(A)$ and $T^*(A^*)$. 

\section{Extension of $\mx$ to arbitrary \dvbs}
\label{sect:dvbsd}

In \cite{MackenzieX:1994}, $\mx$ was used to transfer the standard cotangent
vector bundle structure of $T^*(A^*)\to A^*$ to $T^*(A)\to A^*$ and this
equipped $T^*(A)$ with the structure of a \dvb with core $T^*M$. It then followed from the
local formula for $\mx$ that it was an isomorphism of \dvbs, preserving
$A$ and $A^*$, and reversing the sign on the cores. 

In this paper we assume the general concept of a double
vector bundle \cite{Pradines:FVD} and the basic
process of passing from a \dvb to its duals, as presented in \cite[\S9.2]{Mackenzie:GT}. 

To recall, the \dvb $T(A)$, shown in the first diagram in
Figure~\ref{fig:TA}, has dual \dvb over $A$, denoted $T(A)\duer A$,
identical with the cotangent \dvb $T^*(A)$. The dual \dvb over $TM$,
denoted $T(A)\duer TM$ here and by $T^\sol(A)$ in \cite{Mackenzie:GT},
is canonically isomorphic to $T(A^*)$ under a map $T(A^*)\to T(A)\duer TM$
associated to the tangent pairing $\llangle~,~\rrangle$ and denoted $I$ in
\cite{Mackenzie:GT}. The dual of this \dvb over $A^*$, denoted
$T(A)\duer TM\duer A^*$, is canonically isomorphic to $T^*(A^*)$. 

\begin{figure}[h]
$$
\begin{tikzcd}[ampersand replacement=\&, row sep = 1cm, column sep = 1cm]
T(A) \arrow[r, "T(q)"]
\arrow[d,"",swap]
\& TM \arrow[d,""]
\\
A\arrow[r,"q",swap] \& M,
\end{tikzcd}
\quad
\begin{tikzcd}[ampersand replacement=\&, row sep = 1cm, column sep = 1cm]
T^*(A)\arrow[r, ""]
\arrow[d,"",swap]
\& A^* \arrow[d,"q_*"]
\\
A\arrow[r,"q",swap] \& M,
\end{tikzcd}
\quad
\begin{tikzcd}[ampersand replacement=\&, row sep = 1cm, column sep = 1cm]
T(A^*)\arrow[r, "T(q_*)"]
\arrow[d,"",swap]
\& TM \arrow[d,""]
\\
A^*\arrow[r,"q_*",swap] \& M, 
\end{tikzcd}
\quad
\begin{tikzcd}[ampersand replacement=\&, row sep = 1cm, column sep = 1cm]
T^*(A^*)\arrow[r, ""]
\arrow[d,"",swap]
\& A \arrow[d,""]
\\
A^*\arrow[r,"q_*",swap] \& M.
\end{tikzcd}
$$
\caption{\label{fig:TA}}
\end{figure}

Equation (\ref{eq:cond-v2}) can be written in `mnemonic form' as
\begin{equation}
\label{eq:cond-v3}
\langle T^*(A^*),T(A^*)\rangle + \langle T^*(A),T(A)\rangle =
\llangle T(A^*),T(A)\rrangle.
\end{equation}
This makes clear that condition (\ref{eq:cond-v2}), which defines
the graph of $\mx$, involves the \dvb $T(A)$, its two duals and an
iterated dual. We now apply this process to an arbitrary \dvb.

In Figure~\ref{fig:Dds}, the fourth diagram is the dual of $D\duer B$ over $C^*$;
note that the side structures are the same as for $D\duer A$ but their positions
are interchanged. Figure~\ref{fig:TA} is the case $D = T(A)$, $B = TM$ (and $C = A$). 

We now establish that the following generalization of equation
(\ref{eq:cond-v3}) defines the graph of an isomorphism of \dvbs. The
subscripts in (\ref{eq:3termDBA}) denote the bases of the pairings. 

\begin{equation}
\label{eq:3termDBA}
\langle D\duer B\duer C^*,D\duer B\rangle_{C^*}
 + \langle D\duer A,D\rangle_A 
 = \langle D\duer B,D\rangle_B. 
\end{equation}

We first change some of the notation. We write $d\in D$, $\Phi\in D\duer A$
and $\Psi\in D\duer B$, as in \cite{Mackenzie:GT} and elsewhere. For elements
of the iterated duals, we write $\mA\in D\duer A\duer C^*$ and
$\mB\in D\duer B\duer C^*$. 

\begin{theorem}
\label{thm:D3term}
Let $D$ be a \dvb as in Figure~\ref{fig:Dds}. 
Given $\mB\in D\duer B\duer C^*$ there is a unique $\Phi\in D\duer A$
such that 
\begin{equation}
\label{eq:D3term}
\langle \mB,\Psi\rangle_{C^*} + \langle\Phi,d\rangle_A 
= \langle \Psi,d\rangle_B, 
\end{equation}
for all compatible $\Psi\in D\duer B$ and $d\in D$. 

This condition defines the graph of an isomorphism of \dvbs
$$
\RM{D}{A}
\co D\duer B\duer C^* \to D\duer A, 
$$
which preserves $A$ and $C^*$ and induces $-\id$ on the cores $B^*\to B^*$; that is, 
\begin{equation}
  \label{eq:3termda}
  \langle \mB,\Psi\rangle_{C^*} + \langle\RM{D}{A}(\mB),d\rangle_A 
  = \langle \Psi,d\rangle_B. 
\end{equation}
for all compatible $\Psi\in D\duer B$ and $d\in D$. 
\end{theorem}

The proof follows the same pattern as the proof of Theorem~\ref{thm:MX94}. 

{\sc Proof:}\ 
Consider pairs $(\mB,\Phi)$ in $(D\duer B\duer C^*)\times (D\duer A)$
that project to the same $m\in M$, as shown in Figure~\ref{fig:els}.

\begin{figure}[h]
$$
\begin{tikzcd}[ampersand replacement=\&, row sep = 1cm, column sep = 1cm]
d \arrow[r,mapsto, ""]
\arrow[d,mapsto,"",swap]
\& b_2 \arrow[d,mapsto,""]
\\
a_3\arrow[r,mapsto,"",swap] \& m,
\end{tikzcd}
\quad
\begin{tikzcd}[ampersand replacement=\&, row sep = 1cm, column sep = 1cm]
\Phi \arrow[r,mapsto, ""]
\arrow[d,mapsto,"",swap]
\& \gamma_2 \arrow[d,mapsto,""]
\\
a_2\arrow[r,mapsto,"",swap] \& m,
\end{tikzcd}
\quad
\begin{tikzcd}[ampersand replacement=\&, row sep = 1cm, column sep = 1cm]
\Psi \arrow[r,mapsto, ""]
\arrow[d,mapsto,"",swap]
\& b_1 \arrow[d,mapsto,""]
\\
\gamma_3\arrow[r,mapsto,"",swap] \& m,
\end{tikzcd}
\quad
\begin{tikzcd}[ampersand replacement=\&, row sep = 1cm, column sep = 1cm]
\mB \arrow[r,mapsto, ""]
\arrow[d,mapsto,"",swap]
\& a_1 \arrow[d,mapsto,""]
\\
\gamma_1\arrow[r,mapsto,"",swap] \& m.
\end{tikzcd}
$$
\caption{\label{fig:els}}
\end{figure}

If the elements in Figure~\ref{fig:els} satisfy (\ref{eq:D3term})
then $a_1 = a_2 = a_3$, $b_1 = b_2$ and $\gamma_1 = \gamma_2 = \gamma_3$.
We can therefore take $D$ to be a decomposed \dvb,
$D = A\times_M B\times_M C$.

We write $D\duer A = A\times_M B^*\times_M C^*$,
$D\duer B = C^*\times_M A^*\times_M B$,
and $D\duer B\duer C^* = C^*\times_M B^*\times_M A$. 

Write $\mB = (\gamma_1,\beta_1,a_1)$, $\Psi = (\gamma_3,\alpha_1,b_1)$,
$\Phi = (a_2,\beta_2,\gamma_2)$, and $d = (a_3,c,b_2)$.

We fix $\mB$ and seek to find $a_2$, $\beta_2$ and $\gamma_2$.

For the pairings in (\ref{eq:D3term}) to be defined we must have
$a_2 = a_3$, $b_1 = b_2$ and $\gamma_1 = \gamma_3$.
Equation (\ref{eq:D3term}) then becomes
\begin{equation}
  \label{eq:6pairings}
(\langle\beta_1,b_1\rangle + \langle\alpha_1,a_1\rangle) +
(\langle \beta_2,b_1\rangle + \langle\gamma_2,c\rangle) =
\langle \gamma_1,c\rangle + \langle \alpha_1,a_2\rangle.
\end{equation}
The core elements $\alpha_1$ and $c$ may take any values,
and so too may $b_1 = b_2$. 

Set $c = 0$ and $\alpha_1 = 0$ in (\ref{eq:6pairings}).
We get $\beta_1 + \beta_2 = 0$. Then setting $c = 0$ gives
$a_1 = a_2$. Lastly, $\gamma_1 = \gamma_2$.

So $\Phi = (a_1,-\beta_1,\gamma_1)$ is uniquely determined by $\mB$.
Conversely, $\Phi$ determines $\mB$ uniquely. So there is a diffeomorphism
$$
\RM{D}{A}
\co D\duer B\duer C^* \to D\duer A, 
$$
given in terms of a decomposition by $\RM{D}{A}(\gamma,\beta,a) = (a,-\beta,\gamma)$. 
It is an isomorphism of \dvbs preserving $C^*$ and $A$, and induces
$-\id$ on the cores $B^*\to B^*$.

Define
\begin{equation}
\RM{D}{B} = \RM{D}{A}\duer C^*\co D\duer A\duer C^* \to D\duer B. 
\end{equation}
This is also an isomorphism of \dvbs, but is $+\id$ on the cores
$A^*\to A^*$ and is $-\id$ on the side $B\to B$.

Using $\RM{D}{A}$ we define a pairing of $D\duer A$ with $D\duer B$.
That is, define
\begin{equation}
\label{eq:nsp}
\newpair{\Phi}{\Psi} := \langle (\RM{D}{A})^{-1}(\Phi),\Psi\rangle_{C^*} 
= \langle \Psi,d\rangle_B - \langle\Phi,d\rangle_A,  
\end{equation}
where $d$ is any compatible element of $D$. The second equality follows from 
(\ref{eq:3termda}). 
	
This is the opposite sign convention to that of (\ref{eq:nspAB}) 
and \cite{Mackenzie:GT,Mackenzie:2011}; we use the
$\newpair{}{}$ notation to emphasize this. In the rest of
the paper we will use (\ref{eq:nsp}). 

\begin{remark}
Note that in terms of local decompositions, the sign conventions
  are reversed. For example, the BA sign convention (\ref{eq:nsp}) gives
\begin{equation}
\label{eq:locnsp}  
\newpair{(a,\beta,\gamma)}{(\gamma,\alpha,b)}
= \langle\alpha,a\rangle - \langle\beta,b\rangle. 
\end{equation}
\end{remark}

\begin{remark}
  \label{rmk:ZZ}
In \cite{Mackenzie:GT,Mackenzie:2011}, maps $Z_A$ and $Z_B$ are defined by
\begin{gather*}
Z_A\co D\duer A\to D\duer B\duer C^*,\qquad
\langle Z_A(\Phi), \Psi\rangle_{C^*} = \lhangle \Phi, \Psi\rhangle\\
Z_B\co D\duer B\to D\duer A\duer C^*,\qquad
\langle Z_B(\Psi), \Phi\rangle_{C^*} = \lhangle \Phi, \Psi\rhangle
\end{gather*}
where $\Phi\in D\duer A$, $\Psi\in D\duer B$ and
$\lhangle~,~\rhangle$ is the \nsp in the AB convention.
Thus $\RM{D}{A}$, defined in terms of the BA convention and $Z_A$, defined
in terms of the AB convention, are related by
\begin{equation}
\label{eq:BAAB}
\langle (\RM{D}{A})^{-1}(\Phi),\Psi\rangle_{C^*} 
= -\langle Z_A(\Phi), \Psi\rangle_{C^*},
\end{equation}
and similarly for $Z_B$ and $\RM{D}{B}$. In what follows
we will use the maps $\RM{D}{A}$ and $\RM{D}{B}$. 
\end{remark}

\section{Warps and grids on \dvbs}
\label{sect:warpsandgrids}

We now enter the second part of the paper.

Warps and grids on double and triple vector bundles
were introduced in \cite{Mackenzie:pJihw}, as a tool
to formulate and prove the Jacobi identity of vector
fields diagrammatically. These tools were further
developed in \cite{Flari:thesis} and in \cite{FlariM:2019}.

In this section we briefly present the necessary definitions
and examples, in order to make the link between the nonstandard
duality of \dvbs and warps and grids on \dvbs in 
Section \ref{ssubsect:sqcaps}.
	
The first principle on which the entire apparatus of warps 
and grids lies, is the fact (which can be proven,
quite easily) that two elements $d$ and $d'$ in a \dvb $D$, as
shown in the first and second diagrams in Figure \ref{fig:els-D},
which project to the same $a\in A$ and the same $b\in B$, 
differ by a \emph{unique} core element $c\in C$, see 
\cite[Section 1.1]{FlariM:2019} for more details.

\begin{figure}[h]
  $$
  \begin{tikzcd}[row sep = 1cm, column sep = 1cm]
  d \arrow[r,mapsto, ""]
  \arrow[d,mapsto,"",swap]
  & b \arrow[d,mapsto,""]
  \\
  a\arrow[r,mapsto,"",swap] & m,
  \end{tikzcd}
  \quad
  \begin{tikzcd}[row sep = 1cm, column sep = 1cm]
    d' \arrow[r,mapsto, ""]
    \arrow[d,mapsto,"",swap]
    & b \arrow[d,mapsto,""]
    \\
    a\arrow[r,mapsto,"",swap] & m,
    \end{tikzcd}  \quad
    \begin{tikzcd}[row sep = 1cm, column sep = 1cm]
     d\sub{A} d' \arrow[r,mapsto, ""]
      \arrow[d,mapsto,"",swap]
      & 0^B_m \arrow[d,mapsto,""]
      \\
      a\arrow[r,mapsto,"",swap] & m,
      \end{tikzcd}  
  \quad
  \begin{tikzcd}[row sep = 1cm, column sep = 1cm]
    d\sub{B} d' \arrow[r,mapsto, ""]
     \arrow[d,mapsto,"",swap]
     & b \arrow[d,mapsto,""]
     \\
     0^A_m\arrow[r,mapsto,"",swap] & m.
     \end{tikzcd}  
 $$
  \caption{\label{fig:els-D}}
  \end{figure}
More specifically, whether we take their difference over $A$, or over $B$, the
core element $c$ defined is unique:
 \begin{equation}
    \label{eqn:unique-core-elmt}
  d\sub{A} d' = c\add{B}  \tilo^A_{a},\qquad
  d\sub{B} d'= c\add{A} \tilo^B_{b}.
 \end{equation}
 (We reserve the notation $\Bar{c}$, used in \cite{Mackenzie:GT} and elsewhere,
 for explicit examples in which the core is identified with a bundle that is
 not a subset of $D$.)
 
Now, consider a \dvb $D$ as in the first diagram in 
Figure \ref{fig:threegrids}.

\begin{definition}
  A pair of sections $(\xi,X)$, where $X\in\Ga A$ and $\xi\in \Ga_BD$,
  form a \emph{linear section} of $D$ if $\xi$ is a morphism of vector bundles over $X$. 

A \emph{grid on $D$} is a pair of linear sections $(\xi, X)$ and $(\eta,Y)$
as shown in the first diagram in Figure~\ref{fig:threegrids}.
\end{definition}  

\begin{figure}[h]
  $$
  \begin{tikzcd}[row sep=1.5cm, column sep = 1.5cm]
    D \arrow[r]\arrow[r,<-, bend left = 25, "\xi"]
    \arrow[d]\arrow[d,<-,bend right = 25,swap, "\eta"] 
    &B \arrow[d]\arrow[d,<-, bend left = 25, "Y"]\\
    A\arrow[r]\arrow[r,<-, bend right = 25,swap, "X"]
    &M.
    \end{tikzcd}
          \quad
          \begin{tikzcd}[row sep=1.5cm, column sep = 1.5cm]
            T^2M \arrow[r, "T(p)",swap]\arrow[r,<-, bend left = 25,"T(Y)"]
            \arrow[d,"p_{TM}"]\arrow[d,<-, bend right = 25,swap, "\Tilde{X}"] &TM
            \arrow[d,"p",swap]\arrow[d,<-, bend left = 25, "X"]\\
            TM\arrow[r,"p"]\arrow[r,<-, bend right = 25,swap, "Y"]
            &M.
            \end{tikzcd} 
     \quad
     \begin{tikzcd}[row sep=1.5cm, column sep = 1.5cm]
      T(A) \arrow[r,"T(q)",swap]\arrow[r,<-, bend left = 25,
      "T(\mu)"]\arrow[d,"p_A"]\arrow[d,<-, bend right = 25,swap, "Z^{\hlift}"] &TM
      \arrow[d,"p",swap]\arrow[d,<-, bend left = 25, "Z"]\\
      A\arrow[r,"q"]\arrow[r,<-, bend right = 25,swap, "\mu"]
      &M.
      \end{tikzcd}
            $$
  \caption{\label{fig:threegrids}}
  \end{figure}

For each $m\in M$, $\xi(Y(m))$ and $\eta(X(m))$ have the same outline. 
Hence they determine a unique element of the core $C$ and, as $m$ varies, 
a section of $C$ which we denote $\warp(\xi, \eta)$. Precisely,
\begin{gather}
  \label{eqn:warp of xi and eta}
  \begin{split}
\xi(Y(m))\sub{A} \eta(X(m)) = \warp(\xi,\eta)(m)\add{B} \tilo^A_{X(m)},\\
\xi(Y(m))\sub{B} \eta(X(m))= \warp(\xi,\eta)(m)\add{A} \tilo^B_{Y(m)}.
  \end{split}
\end{gather}

\begin{definition}
\label{df:warp}
The \emph{warp} of the grid on $D$ consisting of
$(\xi, X)$ and $(\eta,Y)$ is $\warp(\xi,\eta)\in\Ga C$. 
\end{definition}

\begin{remark}\label{rmk:sign-warp}
  A comment worth making, which we will use in the proof
  of Theorem \ref{thm:warpsqcap}, is that the 
  warp $\warp(\xi,\eta)$ changes sign if $\xi$ and $\eta$ are 
  interchanged. Our convention gives the positive sign to the 
  counterclockwise composition $\xi\circ Y$. 
  When we take the clockwise composition, the sign of the
  warp changes.    
\end{remark}

\begin{example}\label{ex:T2M}
 In our first example, we take $D = T^2M = T(TM)$, the 
  \emph{double tangent bundle} of a manifold $M$. This is a 
  special case of $T(A)$, with $A = TM$. In this 
  \dvb there are three copies of $TM$; two are
  the side bundles and the third copy is the core.

  Take the following grid on $T^2M$ as shown in
  the second diagram in Figure \ref{fig:threegrids}, 
  comprised of the linear sections $(T(Y),Y)$
  and $(\Tilde{X},X)$, where $X$ and $Y$ are vector fields 
  on $M$. Here $\Tilde{X}$ denotes the complete lift of $X$ to a vector field on $TM$. 
  The complete lift, or tangent lift, $\Tilde{X}$ is $J_M\circ T(X)$ where
  $J_M: T^2M\to T^2M$ is the canonical 
  involution which interchanges the two bundle structures on $T^2M$,
  see \cite[Section 9.6]{Mackenzie:GT} for more details
  on the map $J_M$.
It can be proved that the warp of this grid is 
the Lie bracket $[X,Y]\in \vfk{M}$, see 
\cite[Section 1.1]{FlariM:2019},
\begin{equation}
  \label{eq:AMRi}
  T(Y)(X(m)) \sub{p_{TM}} \Tilde{X}(Y(m)) = ([X,Y])^\uparrow(Y(m)) = 
  [X,Y](m)\add{T(p)} \tilo_{Y(m)},
\end{equation}
where the final zero is for the standard tangent bundle structure,
and the uparrow denotes the vertical lift to $TM$ of the vector 
$[X,Y](m)$ to $Y(m)$. 
\end{example}
  
\begin{example}\label{ex:conn}
Consider the tangent \dvb $T(A)$ of a \vb $A\to M$, and
let $\nabla$ be a connection in $A$.
Take a vector field $Z$ on $M$, and take its 
horizontal lift $Z^{\hlift}$ to $A$ induced by $\nabla$, defined by
\begin{equation}\label{eqn:hlift-linear-pullbacks}
  Z^\hlift(\ell_\ph) = \ell_{\nabla^{(*)}_Z(\ph)},
  \qquad Z^\hlift(f\circ q) = Z(f)\circ q,   
\end{equation} 
where $\nabla^{(*)}$ is the connection in $A^*$ dual to $\nabla$.
Here we have used the fact that to define a vector field on $A$, 
it is sufficient to define its effect on pullback functions 
$f\circ q$ for $f\in C^{\infty}(M)$, and on linear
functions $\ell_{\ph}$ for $\ph \in\Ga A^*$, defined by
$\ell_\ph(a) = \langle\ph(q(a)),a\rangle$, for $a\in A$. 
  
Since $Z^\hlift$ maps linear functions to linear functions
and pullbacks to pullbacks, it is a linear vector field; it clearly
projects to $Z$. The word `horizontal' here has its standard meaning
in connection theory, and does not refer to the structures in $T(A)$.
   
Now take any $\mu\in\Ga A$ and form the grid shown in
the third diagram in Figure~\ref{fig:threegrids}.
The warp of the grid is $\nabla_Z\mu$; that is, for $m\in M$, 
  \begin{equation}
  \label{eq:conn}
  T(\mu)(Z(m)) - Z^\hlift(\mu(m)) = ((\nabla_Z\mu)(m))^\uparrow(\mu(m)),
  \end{equation}
  where the right hand side is the vertical lift of
  $(\nabla_Z\mu)(m)\in A_m$ to $T_{\mu(m)}A$, see \cite[Example~1]{FlariM:2019}  
  for more details.

  More generally let $(\xi, X)$ be any linear vector field on $A$
  and write $D:\Gamma A\to \Gamma A$ for the corresponding differential
  operator (\cite[\S3.4]{Mackenzie:GT}). Then for any $\mu\in\Gamma A$
  the grid formed by $(\xi, X)$ and $(T(\mu),\mu)$ has warp $D(\mu)\in\Gamma A$. 
  \end{example}
	
\section{Warps and duality}
\label{sect:app}
We now prove Theorem \ref{thm:warpsqcap}, an alternative 
formula for the warp which relies on the duality of \dvbs. In
subsections \ref{subsect:gwT2M} and \ref{subsect:gwTA} we
verify Theorem \ref{thm:warpsqcap} directly for the grids 
in the second and third diagrams in Figure~\ref{fig:threegrids}.

\subsection{Squarecap sections and their pairings}
\label{ssubsect:sqcaps}
 
The following material on linear sections and duality can be found in more 
detail in \cite{Mackenzie:2011}.

\begin{figure}[h]
  \begin{tikzcd}[row sep=1cm, column sep = 1cm]
D \arrow[r]\arrow[r,<-, bend left = 25, "\xi"]\arrow[d]\arrow[d,<-,
bend right = 25,swap, "\eta"] &B \arrow[d]\arrow[d,<-, bend left =
25, "Y"]\\
A\arrow[r]\arrow[r,<-, bend right = 25,swap, "X"]
&M,
\end{tikzcd}
\quad
\begin{tikzcd}[row sep=1cm, column sep = 1cm]
D\duer A\duer C^* \arrow[r]\arrow[r,<-, bend left = 15,
"\eta^{\sqcap}"]\arrow[d,xshift = 3mm]
&C^*
\arrow[d]
\\
\phantom{\duer C^*}B\arrow[r]\arrow[r,<-, bend right = 15,swap, "Y"]
&M,
\end{tikzcd}\quad
\begin{tikzcd}[row sep=1cm, column sep = 1cm]
D\duer B\duer C^* \arrow[r]
\arrow[d]\arrow[d,<-, bend right = 25,swap, "\xi^{\sqcap}"] &A
\arrow[d]\arrow[d,<-, bend left = 25, "X"]\\
C^*\arrow[r]
&M.
\end{tikzcd}
\caption{\label{diag:D-and-iterated-duals}}
\end{figure}
A linear section $(\eta,Y)$ of $D$ as in the first diagram
in Figure~\ref{diag:D-and-iterated-duals} induces a linear map
$\ell_\eta:D\duer A\to\R$. This map is automatically linear 
with respect to $D\duer A\to A$ but is also linear with respect to 
the other \vb structure, $D\duer A\to C^*$ (\cite[3.1]{Mackenzie:2011}). 
It therefore induces a section of $D\duer A\duer C^*\to C^*$, 
the corresponding \emph{squarecap section}, which we 
denote by $\eta^\sqcap$. This $\eta^\sqcap$ is again a linear 
section over $Y\in\Ga B$.

In more detail, given $\kappa\in C^*$ we define 
$\eta^{\sqcap}(\kappa)\in D\duer A\duer C^*\Big\rvert_{\kappa}$ by defining
its pairing with any $\Phi\in D\duer A\Big\rvert_{\kappa}$ to be
\begin{equation}\label{definition_of_eta_square_cap}
\langle \eta^{\sqcap}(\kappa),\Phi\rangle_{C^*} :=\ell_{\eta}(\Phi) =
\langle\Phi, \eta(\gamma^A_A(\Phi))\rangle_A.
\end{equation}
(We use notations such as $\Big\rvert_{\kappa}$ on \dvbs when the symbol
for the base point makes clear which structure is meant.)
Further, the restriction of $\ell_{\eta}$ to the core of $D\duer A$ is
$\ell_Y:B^*\rightarrow \R$ and $\ell_{\eta^{\sqcap}} = \ell_{\eta}$,
where the first $\ell$ refers to $D\duer B\duer C^*\to C^*$. 

Therefore, we see that there exists a one-to-one correspondence 
between linear sections $(\eta, Y)$ of $D\rightarrow A$ and 
linear sections $(\eta^{\sqcap},Y)$ of 
$D\duer A\duer C^*\rightarrow C^*$.

Similarly, there exists a one-to-one correspondence between linear
sections $(\xi, X)$ of $D\rightarrow B$ and linear sections $(\xi^{\sqcap},X)$
of $D\duer B\duer C^*\rightarrow C^*$, given by 
\begin{equation}\label{definition_of_xi_square_cap}
\langle \xi^{\sqcap}(\kappa),\Psi\rangle_{C^*} :=\ell_{\xi}(\Psi) =
\langle\Psi, \xi(\gamma^B_B(\Psi))\rangle_B, 
\end{equation}
where $\kappa\in C^*$ and $\Psi\in D\duer B\Big\rvert_{\kappa}$. 

\begin{example}
  Let $D = T^2M$ be the double tangent bundle of a manifold $M$
  and consider the grid shown in the second diagram of
  Figure~\ref{fig:threegrids}. For $Y$ a vector field on $M$,
  $T(Y)$ is a linear section of $T(p)\co T^2M\to TM$ over $Y$.
  We show in subsection \ref{subsect:gwT2M} below that $(T(Y)^{\sqcap},Y)$
  can be identified with $(d\ell_Y,Y)$, a linear $1$-form on $T^*M$.

  On the other hand the complete lift $\Tilde{X}$ of a vector
  field $X$ is a linear section of $p_{TM}\co T^2M\to TM$. We
  show in subsection \ref{subsect:gwT2M} that $(\Tilde{X}^\sqcap,X)$
  can be identified with minus the Hamiltonian vector field on
  $T^*M$ defined by $\ell_X$. 
\end{example}

We now transform the pairing $\newpair{}{}$ of $D\duer A$ and $D\duer B$
(see (\ref{eq:nsp})) into a pairing $\lpair~,~\rpair$ between
$D\duer B\duer C^*\to C^*$ and $D\duer A\duer C^*\to C^*$. Given a 
grid in $D$, applying this pairing to the sections defined
in (\ref{definition_of_eta_square_cap}) and
(\ref{definition_of_xi_square_cap}) will give an alternative
formula for the warp. 

Since $D\duer A\duer C^*$ and $D\duer A$ are dual \vbs over $C^*$,
we have the standard pairing between them. We will use this pairing
and the map $\RM{D}{A}$ as in (\ref{eq:3termda}) to define a pairing
between $D\duer A\duer C^*$ and $D\duer B\duer C^*$ over $C^*$.

Take elements $\mB\in D\duer B\duer C^*$ and $\mA\in D\duer A\duer C^*$ with outlines
$(\mB;\kappa,X;m)$ and $(\mA;Y,\kappa;m).$ Using the map
$\RM{D}{A}: D\duer B\duer C^*\rightarrow D\duer A$ we define
\begin{equation}\label{pairing_via_Z_A}
\lpair \mB,\mA\rpair = 
\langle {\mA},\RM{D}{A}(\mB)\rangle_{C^*}. 
\end{equation}

The proof of the next result will take us to the end of the subsection.

\begin{theorem}
	\label{thm:warpsqcap}
	Let $(\xi, X)$ and $(\eta,Y)$ be linear sections forming a grid on a \dvb $D$. Then
	\begin{equation*}
	\lpair\xi^{\sqcap},\eta^{\sqcap}\rpair = \ell_{-\warp(\xi,\eta)}. 
	\end{equation*}
\end{theorem}

\textsc{Proof:}
	The following outlines may help us to keep track of the various calculations,
	\begin{eqnarray*}
		\begin{tikzcd}[row sep=1.5cm, column sep = 1cm]
		  D\duer A\ni \RM{D}{A}(\xi^{\sqcap}(\kappa))
     \arrow[r,mapsto," "]\arrow[d,mapsto,swap," "]
     &\kappa \arrow[d,mapsto," "]\\
		X(m)\arrow[r,mapsto,swap," "]
			&m, 
		\end{tikzcd}&\quad&
\begin{tikzcd}[row sep=1.5cm, column sep = 1cm]
  D\duer A\duer C^*\ni \eta^{\sqcap}(\kappa) 
  \arrow[r,mapsto," "]\arrow[d,mapsto, xshift = 10mm, swap," "] &
  \kappa\arrow[d,mapsto," "]\\
  \phantom{D\duer A\duer C^*\ni }Y(m)\arrow[r,mapsto,swap," "]
  &m. 
\end{tikzcd}
\end{eqnarray*}

        We can now begin the calculation. Start with (\ref{pairing_via_Z_A}),
        with $\xi^{\sqcap}(\kappa)$ in place of $\mB$,
        and $\eta^{\sqcap}(\kappa)$ in place of $\mA$,
	\begin{eqnarray}
          \nonumber
		\lpair \xi^{\sqcap}(\kappa),\eta^{\sqcap}(\kappa)\rpair &=& 
		\langle \eta^{\sqcap}(\kappa),\RM{D}{A}(\xi^{\sqcap}(\kappa))\rangle_{C^*}\\
		&\equals\limits^{(\ref{definition_of_eta_square_cap})}&
		\ell_{\eta}(\RM{D}{A}(\xi^{\sqcap}(\kappa))) \nonumber\\
		&= & \langle \RM{D}{A}(\xi^{\sqcap}(\kappa)),\eta(X(m))\rangle_A.
                          \label{eq:prevcalc}
	\end{eqnarray}

  We now use (\ref{eq:3termda}), for appropriate
  $\Psi\in D\duer B$, that is with outline $(\Psi;\kappa,Y(m);m)$
  as it follows from the proof of Theorem~\ref{thm:D3term}.
  Then (\ref{eq:prevcalc}) becomes
  $$
  \langle \RM{D}{A}(\xi^{\sqcap}(\kappa)),\eta(X(m))\rangle_A = 
  \langle \Psi, \eta(X(m))\rangle_B - 
  \langle \xi^{\sqcap}(\kappa),\Psi\rangle_{C^*},
  $$
and using (\ref{definition_of_xi_square_cap}) for the 
second term in the last equation, we obtain,
  $$
  \langle \RM{D}{A}(\xi^{\sqcap}(\kappa)),\eta(X(m))\rangle_A = 
  \langle \Psi, \eta(X(m))\rangle_B - 
  \langle \Psi,\xi(Y(m))\rangle_B,
  $$
  since $\gamma^B_B(\Psi) = Y(m)$.
Now recall Remark \ref{rmk:sign-warp}, from which it follows that 
$$
\eta(X(m))\sub{B}\xi(Y(m)) = -\warp(\xi,\eta)(m)\add{A}\tilo^B_{Y(m)},
$$
therefore,
\begin{eqnarray}\label{comparison}
	\lpair \xi^{\sqcap}(\kappa),\eta^{\sqcap}(\kappa)\rpair &=& \nonumber
  \langle \RM{D}{A}(\xi^{\sqcap}(\kappa)),\eta(X(m))\rangle_A \\\nonumber
  &=&  \langle \Psi, \eta(X(m))\rangle_B - \langle \Psi,\xi(Y(m))\rangle_B \nonumber\\
  &=&
  \langle \Psi, \eta(X(m))\sub{B}\xi(Y(m))\rangle_B
  \nonumber\\
  & = &
  \langle\Psi, -\warp(\xi,\eta)(m)\add{A}\tilo^B_{Y(m)}\rangle_B.
	\end{eqnarray}

	We now need the formula for the unfamiliar projection
  $\gamma^B_{C^*}:D\duer B\to C^*$, as in the third 
  diagram in Figure~\ref{fig:Dds}. From \cite[\S9.2]{Mackenzie:GT}, this is
	\begin{equation}\label{eq:unfproj}
	\langle \gamma^B_{C^*}(\Psi), c\rangle_M =
	\langle\Psi, \tilo^B_b \add{A} c\rangle,
	\end{equation}
	where $c\in C_m,\ \Psi:(\gamma^B_{C^*})^{-1}(b)\to\R$ and $b\in B_m.$ 
	We can therefore rewrite the last expression of (\ref{comparison}) as
	\begin{equation*}
    \langle\Psi, -\warp(\xi,\eta)(m)\add{A}\tilo^B_{Y(m)}\rangle_B = \langle
	\kappa,-\warp(\xi,\eta)(m)\rangle = -\ell_{\warp(\xi,\eta)}(\kappa), 
	\end{equation*}
	since $\gamma^B_{C^*}(\Psi) = \kappa$. 
	In other words, we have shown that, for $\kappa\in C^*$,
	$\lpair\xi^{\sqcap},\eta^{\sqcap}\rpair(\kappa) = \ell_{-\warp(\xi,\eta)}(\kappa).$ 
	
	This completes the proof of Theorem~\ref{thm:warpsqcap}.

Theorem~\ref{thm:warpsqcap} has recently been found independently by Meinrenken and
Pike \cite[Prop.~3.6]{MeinrenkenP:arx}, using a different approach.

\subsection{Example with $T^2M$}
\label{subsect:gwT2M}

Consider the \dvb $T^2M$ and the grid consisting of $(\Tilde{X}, X)$
and $(T(Y),Y)$, as in the second diagram in Figure~\ref{fig:threegrids}. 
We verify Theorem \ref{thm:warpsqcap} in this case.

\begin{figure}[h]
  $$
  \begin{tikzcd}[row sep=1cm, column sep = 1cm]
    T^2M \arrow[r, "T(p)"] \arrow[d,"p_{TM}"] & TM \arrow[d,"p"] \\
    TM \arrow[r,"p"] &M,
    \end{tikzcd}
      \quad
      \begin{tikzcd}[row sep=1cm, column sep = 1cm]
        T^\bullet TM \arrow[r, "T(p)_\bullet"] \arrow[d,"r_\bullet"] & TM \arrow[d,"p"] \\
     T^*M\arrow[r,"c_M"] &M, 
     \end{tikzcd}
     \quad
     \begin{tikzcd}[row sep=1cm, column sep = 1cm]
      T^*(TM) \arrow[r, "r"] \arrow[d,"c_{TM}"] &T^*M \arrow[d,"c_M"] \\
      TM\arrow[r,"p"] &M.
      \end{tikzcd}
        $$
  \caption{\label{fig:T2Mduals}}
  \end{figure}

First, what are the corresponding squarecap sections 
$({\Tilde{X}}^{\sqcap},X)$ and $(T(Y)^{\sqcap},Y)$? 

We first calculate $T(Y)^{\sqcap}$ using $\ell_{T(Y)}: T^{\bullet}(TM)\to \R$,
where the \dvb $T^{\bullet}(TM)$ (second diagram in Figure~\ref{fig:T2Mduals}),
is the prolongation dual of $T^2M$, the special case $T^{\bullet}(A)$ for 
$A = TM$ as mentioned in Section \ref{sect:dvbsd}. 

As noted in the previous subsection, the function $\ell_{T(Y)}$ is 
linear  with respect to both side bundles of $T^{\bullet}(TM)$. 
Specifically, its linearity with respect to $T^*M$ defines a linear 
section of the \vb $T^{\bullet}(TM)\duer T^*M\to T^*M$, which is 
the dual over $T^*M$ of $T^{\bullet}(TM)\to T^*M$. 

We now use the map $I: T(T^*M)\to T^{\bullet}(TM)$ to canonically
identify $T^{\bullet}(TM)$ with $T(T^*M)$, as mentioned
in Section \ref{sect:dvbsd}. This map $I$, for $A = TM$, 
is defined by 
\begin{equation}
  \label{eq:9.3.2}
  \langle I(\mathcal{X}), \eta\rangle_{TM} 
  = \llangle \mathcal{X}, \eta\rrangle
  \end{equation}
where $\mathcal{X}\in T(T^*M)$ and $\eta\in T^2M$, and again, 
 $\llangle~,~\rrangle$ is the tangent prolongation of the pairing 
 of $TM$ with $T^*M$, as in the proof of Theorem \ref{thm:MX94}.

The map $I$ also allows 
us to simplify the iterated dual $T^{\bullet}(TM)\duer T^*M\to T^*M$, 
which can be identified canonically with the dual over $T^*M$ of 
the \vb $T(T^*M)\rightarrow T^*M$, namely, $T^*(T^*M)\to T^*M$.
Therefore, we are now seeking a linear section $\mathfrak{Y}$ of
$T^*(T^*M)\to T^*M$. This section $\mathfrak{Y}$ is defined by the 
function $\ell_{T(Y)}\circ I: T(T^*M)\rightarrow \R$, 
specifically, by its linearity with respect to $T^*M$.

Consider now $\mathfrak{Y}(\ph)\in T^*(T^*M)$ for $\ph\in T^*M$, and pair 
it with a $\xi\in T(T^*M)$ with outline $(\xi;\ph,x;m)$, where $x\in TM$. Using
(\ref{definition_of_eta_square_cap}), 
\begin{equation*}
\langle \mathfrak{Y}(\ph),\xi\rangle_{T^*M} = 
(\ell_{T(Y)}\circ I)(\xi) = 
\langle I(\xi),T(Y)(x)\rangle_{TM}.
\end{equation*}

We now need the following result from \cite{MackenzieX:1994}; see also 
\cite[3.4.6]{Mackenzie:GT}. It is valid for an arbitrary \vb~$A$. 
  
\begin{proposition}                    
	\label{prop:llangle}
	Given $(\xi;\mu(m),x;m)\in T(A)$ and $(\mathfrak{X};\ph_m,x;m)\in T(A^*),$ 
	let $\mu \in\Gamma(A)$ and $\ph\in\Gamma(A^*)$ be any sections taking 
	the values $\mu(m)$ and $\ph_m$ at $m.$ Then
	\begin{equation}
	\label{eq:llangle}
	\llangle\mathfrak{X},\xi\rrangle = \mathfrak{X}(\ell_\mu) + \xi(\ell_\ph)
	- x(\langle\ph,\mu\rangle).
	\end{equation}
\end{proposition}

Using Proposition \ref{prop:llangle} and (\ref{eq:9.3.2}), it follows that
\begin{equation}\label{eq:def-of-T(Y)}
\langle I(\xi),T(Y)(x)\rangle_{TM} = 
\llangle \xi, T(Y)(x)\rrangle  
= \xi(\ell_Y) =\langle d\ell_Y(\ph), \xi\rangle.
\end{equation}
This is true for any such $\xi\in T(T^*M)$, so it follows that
\begin{equation*}
\mathfrak{Y}(\ph) = (d\ell_Y)(\ph),
\end{equation*} 
and the linear section in question, $(T(Y)^{\sqcap},Y)$,
 can be identified with $(d\ell_Y,Y)$.

\medskip

Next consider the linear section $(\Tilde{X},X)$ and the third \dvb in 
Figure~\ref{fig:T2Mduals}. 
As before, the function
\begin{equation*}
\ell_{\Tilde{X}}: T^*(TM)\rightarrow \R.
\end{equation*}
is linear with respect to both $TM$ and $T^*M$. 
Using the linearity over $T^*M$, we obtain a section of the dual over $T^*M$
of the \vb $T^*(TM)\to T^*M$; that is, of $T^*(TM)\duer T^*M\to T^*M$.

Again, this is not easy to work with, and in this case we need to use 
the map $R:T^*(A^*)\rightarrow T^*(A)$ as in Theorem \ref{thm:MX94}, for 
the case $A = TM$. It follows that
\begin{equation*}
\ell_{\Tilde{X}}\circ R: T^*(T^*M)\to \R
\end{equation*}
is also linear with respect to $T^*M$,
and defines a section $\mathfrak{X}$ of the dual of the \vb
$T^*(T^*M)\to T^*M$; that is, of the tangent bundle $T(T^*M)\to T^*M$.

Then for $\ph \in T^*M$, and any $\mathfrak{F}\in T^*(T^*M)$ with outline
$(\mathfrak{F};x,\ph;m)$, with $x\in TM$, using
(\ref{definition_of_eta_square_cap}),
\begin{equation*}
\langle \goX(\ph),\mathfrak{F}\rangle_{T^*M} = (\ell_{\Tilde{X}}\circ
R)(\mathfrak{F}) = \ell_{\Tilde{X}}(R(\mathfrak{F})) = \langle
R(\mathfrak{F}),\Tilde{X}(x)\rangle_{TM}.
\end{equation*}

At this point, we use the commutative diagram (\ref{eq:RJIdnu}), in which
each map is an isomorphism of \dvbs and $(d\nu)^{\sharp}$ is the map associated to
the canonical symplectic structure $d\nu$ on $T^*M$; see \cite{MackenzieX:1994}
or \cite[p. 442]{Mackenzie:GT}. 
\begin{equation}
\label{eq:RJIdnu}  
\begin{tikzcd}[row sep=1.25cm, column sep = 1.25cm]
T^*(T^*M) \arrow[r,"R"]\arrow[d, swap,"(d\nu)^{\sharp}"] 
&T^*(TM) \arrow[d,<-,"J^*"]\\
T(T^*M)\arrow[r,swap,"I"]
&T^{\bullet}(TM). 
\end{tikzcd}
\end{equation}
Using $R = J^*\circ I \circ (d\nu)^{\sharp}$, we have
\begin{multline*}
\langle R(\mathfrak{F}),\Tilde{X}(x)\rangle_{TM} = 
\langle J^*(I ((d\nu)^{\sharp}(\mathfrak{F}))), 
\Tilde{X}(x)\rangle_{TM} \\
= \langle I ((d\nu)^{\sharp}(\mathfrak{F})),  J(\Tilde{X}(x))\rangle_{TM} = 
\langle I ((d\nu)^{\sharp}(\mathfrak{F})),  T(X)(x)\rangle_{TM}.
\end{multline*}
As before, using (\ref{eq:def-of-T(Y)}),
\begin{equation*}
\langle I ((d\nu)^{\sharp}(\mathfrak{F})),  T(X)(x)\rangle_{TM}
= (d\nu)^{\sharp}(\mathfrak{F})(\ell_X) =
\langle d\ell_X,(d\nu)^{\sharp}(\mathfrak{F}) \rangle = 
-\langle (d\nu)^{\sharp}(d\ell_X),\mathfrak{F}\rangle,
\end{equation*}
so we see that $\mathfrak{X} = -(d\nu)^{\sharp}(d\ell_X)$; that is, it is minus 
the Hamiltonian vector field for the function $\ell_X$. Denote it by
$H_{\ell_X}$. Finally,
\begin{equation*}
\lpair T(Y)^{\sqcap}, \Tilde{X}^{\sqcap}\rpair_{T^*M} = \langle d\ell_Y,
-H_{\ell_X}\rangle = -\ell_{[X,Y]}, 
\end{equation*}
and the warp of the grid is $[X,Y]$ as in Example~\ref{ex:T2M}.

\subsection{Example with $T(A)$}
\label{subsect:gwTA}

In the case of the grid on $T(A)$ as in 
the third diagram in Figure~\ref{fig:threegrids}, 
what are the corresponding sections $T(\mu)^{\sqcap}$ and $(X^{\hlift})^{\sqcap}$~?
Just as in the case of $T^2M$, $(T(\mu)^{\sqcap},\mu)$ can be identified with
$(d\ell_{\mu},\mu)$.

For $(X^{\hlift},X)$ a more elaborate calculation is needed. 
Again, we use $\ell_{X^{\hlift}}\circ R: T^*(A^*)\rightarrow \R$, and its
linearity with respect to $A^*$. This will define a section of the dual of
$T^*(A^*)\rightarrow A^*$, that is, of $T(A^*)\rightarrow A^*$. Denote this
vector field by $\Phi$.

Given $\kappa\in A^*$, pair $\Phi(\kappa)\in T(A^*)$ with any $\Psi\in T^*(A^*)$
which has outline $(\Psi;\kappa,a;m)$. 

Using equation (\ref{eq:cond-v2}), it follows that, for suitable 
$\mathcal{X}\in T(A^*)$, 
\begin{equation}
\label{eq:hence}
\langle \Phi(\kappa),\Psi\rangle_{A^*} = 
(\ell_{X^{\hlift}}\circ R)(\Psi) = 
\langle R(\Psi), X^{\hlift}(a)\rangle_A = 
\llangle \mathcal{X}, X^{\hlift}(a)\rrangle - \langle \Psi,
\mathcal{X}\rangle_{A^*}.
\end{equation}
The outlines for the elements involved are:
\begin{equation*}
\begin{tikzcd}[row sep=1cm, column sep = 1cm]
T^*(A)\ni R(\Psi) \arrow[r,mapsto," "]\arrow[d,mapsto, swap," ",xshift
= 4.5mm] &a\in A \arrow[d,mapsto," ",xshift = -3.75mm]\\
A^*\ni \kappa\arrow[r,mapsto,swap," "]
&m,\phantom{\in A}
\end{tikzcd}\quad
\begin{tikzcd}[row sep=1cm, column sep = 1cm]
T(A)\ni X^{\hlift}(a)\arrow[r,mapsto," "]\arrow[d,mapsto, swap," ",xshift =
4.5mm] &X(m)\in TM \arrow[d,mapsto," ",xshift = -5.75mm]\\
A\ni a\arrow[r,mapsto,swap," "]
&m,\phantom{\in TM} 
\end{tikzcd}
\end{equation*}
\begin{equation*}
\begin{tikzcd}[row sep=1cm, column sep = 1cm]
T^*(A^*)\ni \Psi \arrow[r,mapsto," "]\arrow[d,mapsto, swap," ",xshift =
6.5mm] &a\in A \arrow[d,mapsto," ",xshift = -3.75mm]\\
A^*\ni \kappa\arrow[r,mapsto,swap," "]
&m,\phantom{\in A}
\end{tikzcd}\quad
\begin{tikzcd}[row sep=1cm, column sep = 1cm]
T(A^*)\ni\mathcal{X} \arrow[r,mapsto," "]\arrow[d,mapsto, swap," ",xshift =
5mm] &X(m)\in TM \arrow[d,mapsto," ",xshift = -5.75mm]\\
A^*\ni \kappa\arrow[r,mapsto,swap," "]
&m.\phantom{\in TM} 
\end{tikzcd}
\end{equation*}
Now use Proposition \ref{prop:llangle} for the first term of (\ref{eq:hence}). 
Choose a $\ph\in \Gamma A^*$ with $\ph(m) = \kappa$, and a $\mu\in \Gamma A$
with $\mu(m) = a$. We can also make the following choice; linear vector fields
of a \vb $A$ are in bijective correspondence with linear vector fields on its
dual bundle $A^*$ (see \cite[3.4.5]{Mackenzie:GT}). Therefore, take
$\mathcal{X}$ to be $X^{\hlift_*}(\ph(m))$, where $X^{\hlift_*}$ is the corresponding
linear vector field to $X^{\hlift}$. Then we can write
\begin{multline}\label{eq:ex-conn-calculations}
\llangle X^{\hlift_*}(\ph(m)), X^{\hlift}(a)\rrangle - \langle \Psi,
X^{\hlift_*}(\ph(m))\rangle_{A^*} \\=
X^{\hlift_*}(\ph(m))(\ell_{\mu}) + X^{\hlift}(\mu(m))(\ell_{\ph})
-X(m)(\langle \ph,\mu\rangle) - \langle \Psi,X^{\hlift_*}(\ph(m))\rangle_{A^*}.
\end{multline}
At this point, recall from (\ref{eqn:hlift-linear-pullbacks}) 
that for $\ph\in \Gamma A^*$, and for $\mu\in \Gamma A$,
\begin{equation*}
X^{\hlift}(\ell_{\ph}) = \ell_{\nabla_X^{\dual}(\ph)}\in
C^{\infty}(A),\qquad
X^{\hlift_*}(\ell_{\mu}) = \ell_{\nabla_X(\mu)}\in C^{\infty}(A^*), 
\end{equation*}
and of course the relation between $\nabla$ and $\nabla^{\dual}$,
\begin{equation*}
\langle \nabla_X^{\dual}(\ph), \mu\rangle = 
X(\langle \ph, \mu\rangle) - \langle\ph, \nabla_X(\mu)\rangle,
\end{equation*}
and the latter equation can be rewritten as
\begin{equation*}
\ell_{\nabla^{(*)}_X(\ph)}\circ \mu = X(\langle \ph, \mu\rangle) -
\ell_{\nabla_X(\mu)}\circ \ph.
\end{equation*}
Returning to (\ref{eq:ex-conn-calculations}),
\begin{multline*}
\llangle X^{\hlift_*}(\ph(m)), X^{\hlift}(a)\rrangle - \langle \Psi,
X^{\hlift_*}(\ph(m))\rangle_{A^*} \\= 
\ell_{\nabla_X(\mu)}(\ph(m))  + \ell_{\nabla^{(*)}_X(\ph)}(\mu(m)) -
X(m)(\langle\ph, \mu\rangle) - \langle \Psi,X^{\hlift_*}(\ph(m))\rangle_{A^*}\\
= -\langle \Psi,X^{\hlift_*}(\ph(m))\rangle_{A^*}.
\end{multline*}

Finally, we have shown that the pairing between
$T(\mu)^{\sqcap}$, which we have shown can be identified with $(d\ell_{\mu},\mu)$, and
$(X^H)^\sqcap$, which can be identified with $-(X^{\hlift_*},X)$ is,
\begin{equation*}
-\langle X^{\hlift_*},d\ell_\mu\rangle =-X^{\hlift_*}(\ell_{\mu}) =
-\ell_{\nabla_X(\mu)}.
\end{equation*}

Again we see that the defining formula for the warp
has been converted into a familiar pairing.

\newcommand{\noopsort}[1]{} \newcommand{\singleletter}[1]{#1} \def\cprime{$'$}
  \def\cprime{$'$}


\begin{thebibliography}{10}

\bibitem{Flari:thesis}
M.~K. Flari.
\newblock {Triple vector bundles in Differential Geometry}.
\newblock Thesis, University of Sheffield, 2018.

\bibitem{FlariM:2019}
M.~K. Flari and K.~Mackenzie.
\newblock Warps, grids and curvature in triple vector bundles.
\newblock {\em Lett. Math. Phys.}, 109(1):135--185, 2019.

\bibitem{KhudaverdyanV:2018}
H.~Khudaverdyan and Th.~Voronov.
\newblock {Thick morphisms, higher Koszul brackets, and
  $L_{\infty}$-algebroids}.
\newblock arXiv:1808.10049.

\bibitem{Mackenzie:pJihw}
K.~Mackenzie.
\newblock Proving the {J}acobi identity the hard way.
\newblock In {\em Geometric methods in physics}, Trends Math., pages 357--366.
  Birkh\"auser/Springer, Basel, 2013.

\bibitem{Mackenzie:1999}
K.~C.~H. Mackenzie.
\newblock On symplectic double groupoids and the duality of {P}oisson
  groupoids.
\newblock {\em Internat. J. Math.}, 10(4):435--456, 1999.

\bibitem{Mackenzie:GT}
K.~C.~H. Mackenzie.
\newblock {\em General theory of {L}ie groupoids and {L}ie algebroids}, volume
  213 of {\em London Mathematical Society Lecture Note Series}.
\newblock Cambridge University Press, Cambridge, 2005.

\bibitem{Mackenzie:2011}
K.~C.~H. Mackenzie.
\newblock Ehresmann doubles and {D}rinfel'd doubles for {L}ie algebroids and
  {L}ie bialgebroids.
\newblock {\em J. Reine Angew. Math.}, 658:193--245, 2011.

\bibitem{MackenzieX:1994}
K.~C.~H. Mackenzie and P.~Xu.
\newblock Lie bialgebroids and {P}oisson groupoids.
\newblock {\em Duke Math. J.}, 73(2):415--452, 1994.

\bibitem{MeinrenkenP:arx}
E.~Meinrenken and J.~Pike.
\newblock {The Weil algebra of a double Lie algebroid}.
\newblock {\sf arXiv:1901.00230}.

\bibitem{Pradines:FVD}
J.~Pradines.
\newblock {\em Fibr{\'e}s vectoriels doubles et calcul des jets non holonomes},
  volume~29 of {\em Esquisses Math\'ematiques [Mathematical Sketches]}.
\newblock Universit\'e d'Amiens, U.E.R. de Math\'ematiques, Amiens, 1977.

\bibitem{Pradines:1988}
J.~Pradines.
\newblock Remarque sur le groupo\"\i de cotangent de {W}einstein-{D}azord.
\newblock {\em C. R. Acad. Sci. Paris S\'er. I Math.}, 306(13):557--560, 1988.

\bibitem{Roytenberg:thesis}
D.~Roytenberg.
\newblock {Courant algebroids, derived brackets and even symplectic
  supermanifolds}.
\newblock Thesis, Univ.~California, Berkeley, 1999. {\sf
  arXiv:math.DG/9910078}.

\bibitem{Tulczyjew:1977ihp}
W.~M. Tulczyjew.
\newblock The {L}egendre transformation.
\newblock {\em Ann. Inst. H. Poincar\'{e} Sect. A (N.S.)}, 27(1):101--114,
  1977.

\bibitem{Tulczyjew:1977}
W.~M. Tulczyjew.
\newblock A symplectic formulation of particle dynamics.
\newblock In {\em Differential geometric methods in mathematical physics
  ({P}roc. {S}ympos., {U}niv. {B}onn, {B}onn, 1975)}, pages 457--463. Lecture
  Notes in Math., Vol. 570. Springer, Berlin, 1977.

\end{thebibliography}
\end{document}